\documentclass{amsart}
\usepackage{mathrsfs}
\usepackage{amsfonts}
\usepackage{txfonts}
\usepackage{amsmath}
\usepackage{amssymb}
\usepackage{amsthm}
\usepackage{graphicx}
\usepackage[toc,page,title,titletoc,header]{appendix}
\usepackage{geometry}
\usepackage{upgreek}
\usepackage[T1]{fontenc}
\usepackage{color}
\usepackage{hyperref}
\usepackage[all]{xy}
\usepackage[french,english]{babel}

\theoremstyle{plain}
\newtheorem{theo}{Theorem}[section]
\newtheorem*{theo*}{Theorem}

\newtheorem*{lem*}{Lemma}

\theoremstyle{definition}
\newtheorem{lem}{Lemma}[section]
\newtheorem{defi}{Definition}[section]

\theoremstyle{remark}

\newtheorem*{ack}{Acknowledgements}
\newtheorem*{rem*}{Remark}
\newtheorem{rem}{Remark}[section]

\newcommand{\R}{\mathbb{R}}

\newcommand{\Z}{\mathbb{Z}}
\newcommand{\C}{\mathbb{C}}

\begin{document}

\title{Generalized periodic orbits in some restricted three-body problems}

\author{Rafael Ortega and Lei Zhao}
%\email{l.zhao@nankai.edu.cn}
%\thanks{Supported by }

\abstract We treat the circular and elliptic restricted three-body problems in inertial frames as periodically forced Kepler problems with additional singularities and explain that in this setting the main result of \cite{BOZ} is applicable. This guarantees the existence of an arbitrary large number of generalized periodic orbits (periodic orbits with possible double collisions regularized) provided the mass ratio of the primaries is small enough. 

\endabstract

\date\today
\maketitle
\section{Introduction}
Consider the system
\begin{equation}\label{eq:1}
\ddot{u}=-\dfrac{u}{|u|^{3}}+ \varepsilon \nabla_{u} U(t, u, \varepsilon), u \in \R^{3}
\end{equation}
where $U$ is a $C^{\infty}$-function defined on $\R \times \mathcal{U} \times [0, \varepsilon_{*}]$ for an open set $\mathcal{U} \subset \R^{3}$ containing the origin $u=0$ and some $\varepsilon_{*}>0$.  The function $U(t, u, \varepsilon)$ is assumed to be $T$-periodic in $t$.  This system models a Kepler problem with an external periodic force. \par

It is a classical problem in the theory of perturbations to look for periodic solutions for small $\varepsilon$. {Traditionally, due to the proper-degeneracy of the Kepler problem, all the collisionless bounded orbits are closed and some non-degeneracy condition has to be imposed on the function $U(t,u,\varepsilon )$. This condition has to be verified for concrete $U(t,u,\varepsilon )$}. In typical situations the periodic orbits are found in the neighbourhood of a prescribed
closed orbit of the unperturbed Kepler problem and in particular they have no collisions.\par 

A result of different nature was obtained in \cite{BOZ}, in which it has been shown that for any smooth function $U$, the equation \eqref{eq:1} has an arbitrarily large number of $T$-periodic solutions if $\varepsilon$
has been accordingly chosen small enough. {As it is usually the case in the study of periodic orbits with singular potentials,} in the framework of \cite{BOZ} the function $U$ was supposed to be globally defined; that is,
\begin{equation}\label{global}
\mathcal{U}=\R^{3}.
\end{equation} Also, the periodic solutions are understood in a generalized sense. They can have collisions but the energy and the direction of motion must be preserved after each bouncing at a time $t_*$ with $u(t_* )=0$.
\par There are good topological reasons for the introduction of generalized solutions. Let $\mathcal{M} $ be the set of $T$-periodic solutions of the Kepler problem ($\varepsilon =0$). After including solutions with collisions it becomes a manifold $$\tilde{\mathcal{M}}=\bigcup_{n=1}^{\infty} \mathcal{M}_n$$ with infinitely many
connected components $\mathcal{M}_n$, where each of them is compact. The periodic solutions of \eqref{eq:1} for $\varepsilon \neq 0$ are obtained as bifurcations from these manifolds. The compactness is essential if we want to guarantee that
this bifurcation is always produced.\par On the other hand, as we willexplain in this note, since our method is perturbative, the condition \eqref{global} is not
so essential. Since the components $\mathcal{M}_n$ converge to the origin as $n\to \infty$, the projection of $\mathcal{M}_n$ on the configuration space will lie inside any neighbourhood $\mathcal{U}$ of the origin for $n$ large enough. In consequence, when \eqref{global}
does not hold,  it is still possible
to obtain bifurcations from $\mathcal{M}_n$ when $n$ is large and the projection of $\mathcal{M}_{n}$ to the configuration space has a neighborhood in which the function $U$ is well-defined. This observation may seem minor but it is of importance if we want to apply this technique to some classical problems in Celestial Mechanics. A typical situation arises in the so-called circular or elliptic restricted three-body problem, where the perturbation $U$ on the asteroid is due to a periodic gravitational force from another massive body. {The function $U$} will then be singular at the position of this massive body which is at some distance of $u=0$. Note that for the circular problem we shall consider the problem in a fixed inertial frame, and shall not go to some rotating frame to reduce the system to an autonomous one.
\par
We shall then apply this to show that an arbitrarily large number of generalized $T$-periodic solutions of a general class of circular or elliptic restricted three body problems with any eccentricity $e \in [0, 1)$  can be found  when the mass ratio between the primaries is taken as an external parameter that we assume to be sufficiently small. These solutions may collide with the big primary.  In fact we will find periodic motions of the small body in a neighbourhood of this big primary. The common year of the primaries will also be a period for the third body and each year this small body will make a large number of revolutions around the big primary. There are many other results on the existence of periodic solutions of the circular or elliptic restricted problem and many families of periodic solutions have been identified. See for instance \cite{CPS, AL, {PYFN}}. A possible novelty of our result is mainly in the absence of additional non-degeneracy conditions and the result is more global in the sense that the continuation is made from manifolds made of periodic orbits of the Kepler problem instead of particular periodic orbits of an approximating system.
It seems reasonable to expect other applications in similar models. For instance, in an elliptic restricted N-body problem where the primaries are assumed to move on an elliptic homographic motion and the infinitesimal body is assumed to stay close to a primary.  
\par
The rest of the paper is organized in three sections. In Section \ref{22} we work with the equation \eqref{eq:1} and go back to the main theorem in \cite{BOZ}. We explain the modifications in the proof allowing to eliminate the condition (\ref{global}). The result in \cite{BOZ} was obtained via the use of Levi-Civita regularization in dimension 2 and Kustaanheimo-Stiefel regularization in dimension 3, together with the application of a theorem of Weinstein \cite{Weinstein}. Generalized periodic solutions can be characterized equivalently by periodic orbits of the regularized system. This was proved in \cite{BOZ} in dimension 2 and we will prove that it is also the case in dimension 3. This proof will require a delicate topological result on the lifting of piecewise smooth paths via the Hopf fibration. For smooth paths this is a direct consequence of general results in the theory of Ehresmann connections, and the modifications required for the piecewise smooth case are explained in Section \ref{33}. Finally, the application to the circular or elliptic restricted three body problem is presented in Section \ref{44}.\par
In this note we will work in the space $u \in \R^{3}$ but the results are easily adapted to the simpler case of the plane $u \in \R^{2}$ or the line $u \in \R$. 

\section{Periodic Solutions via Regularization}\label{22}
Following \cite{BOZ} we say that a continuous and $T$-periodic function $u: \R \to \R^{3}$ is a generalized $T$-periodic solution of \eqref{eq:1} if it satisfies the following conditions:
\begin{enumerate}
\item $\mathcal{Z}=\{t \in \R: u(t)=0\}$ is a discrete set;
\item for any open interval $I \subset \R \setminus \mathcal{Z}$ the function u is in $C^{\infty}(I, \R^{3})$ and satisfies \eqref{eq:1} on this interval;
\item For any $t_{0} \in \mathcal{Z}$ the limits below exist
$$\lim_{t \to t_{0}} \dfrac{u(t)}{|u(t)|}, \qquad \qquad \lim_{t \to t_{0}} \Bigl( \dfrac{1}{2} |\dot{u}(t)|^{2}-\dfrac{1}{|u(t)|} \Bigr).$$
\end{enumerate}

Note that for any classical solution of \eqref{eq:1} tending to a collision at $t_{0}$, the left and right limits $t \to t_{0}^{-}$ and $t \to t_{0}^{+}$ always exist (see \cite{Sperling}). The crucial point in the above definition is that they coincide. 

The above definition extends in an obvious way for the case of sub-harmonic solutions having period $\eta T$ with $\eta \in \Z, \eta \ge 2$.
\par We will prove the existence of generalized $T$-periodic solutions using a regularization technique. We refer to \cite{BDP} and \cite{BOV} for an alternative use of variational techniques. 
Following \cite{ZhaoKS, BOZ} we consider the Kustaanheimo-Stiefel regularization. The skew-field of quaternions is denoted by $\mathbb{H}$. The space of purely imaginary quaternions $\mathbb{IH}=\{z \in \mathbb{H}: \Re(z)=0\}$ is naturally a three-dimensional vector field over the reals. The map 
$$ \Pi:  \mathbb{H} \to \mathbb{IH}, \qquad z \mapsto \bar{z} i z$$ 
plays an important role. 

We set $\mathbb{T}=\R/T \Z$ for the quotient space of the real line of time by the lattice $T \Z$, so for $t \in \R$ we denote by $\bar{t}$ its corresponding quotient. The manifold 
$$\mathcal{M}=\mathbb{H} \times \mathbb{H} \times \mathbb{T} \times \R$$
with points $(z=z_{0}+z_{1} i + z_{2} j + z_{3} k, w=w_{0}+w_{1} i + w_{2} j + w_{3} k, \bar{t}, \tau)$ is endowed with the symplectic two-form
$$\omega=\sum_{l=0}^{3} d z_{l} \wedge d w_{l} + d t \wedge d \tau$$
On the symplectic manifold $(\mathcal{M}, \omega)$ we consider the Hamiltonian function
$$K_{\varepsilon}: \mathcal{M} \to \R, K_{\varepsilon}(z, w, t,\tau)=\dfrac{|w|^{2}}{8}+\tau |z|^{2} -1 + \varepsilon P(t, z, \varepsilon)$$
with $P(t,z,\varepsilon)=|z|^{2} U(t, \bar{z} i z, \varepsilon)$. This is the regularized system of the spatial forced Kepler problem in extended phase space \cite{BOZ}. 

This function is invariant under the Hamiltonian $S^{1}$-action
$$g \ast (z, w, \bar{t}, \tau)=(g z, g w, \bar{t}, \tau),$$
where we realize $S^{1}$ as $\{g \in \C: |g|=1\}$. By standard theory, the corresponding moment map
$$BL(z, w, \bar{t}, \tau)=\bar{z} i w$$
is a first integral of the system. 

Assume now that $X(s)=(z(s), w(s), \bar{t}(s), \tau(s))$ is a solution of the Hamiltonian system $(\mathcal{M}, \omega, K_{\varepsilon})$ which lies in $K_{\varepsilon}^{-1}(0) \cap BL^{-1}(0)$ and $X(s+S)=g \ast X(s), s \in \R$ for some $S >0$ and $g \in S^{1}$. Then $X(s)$ gives a generalized periodic solution of \eqref{eq:1}, as a consequence of Lem 5.1, Rem 5.2 in \cite{BOZ}, with period $\eta T$ where $\eta$ is the degree of the map $\bar{t}(s)$, given via a lift of this mapping as 
$$t(s+S)=t(s)+\eta T.$$

These solutions were produced by a perturbation argument in the symplectically reduced manifold
$$\mathcal{M}_{0}=BL^{-1}(0)/S^{1}.$$
Indeed by applying a result of Weinstein \cite{Weinstein}, we obtained a continuation from the sets $\bar{\Lambda}_{k}$ of periodic orbits of $K_{0}$ which are obtained as quotients from the sets
$$\Lambda_{k}=\{X=(z, w, \bar{t}, \tau) \in K_{0}^{-1}(0) \cap BL^{-1}(0): \tau=\tau_{k}\}$$
with $\tau_{k}=\Bigl(\dfrac{\sqrt{2} k \pi}{T}\Bigr)^{\frac{2}{3}}, \quad k=1,2,\cdots$

In principle the set 
$$E_{k}=\{\bar{z} i z: X=(z, w, \bar{t}, \tau) \in \Lambda_{k}\}$$
can occupy any region in $\mathbb{IH}$. For this reason, it was assumed in \cite{BOZ} that (\ref{global}) holds. Nevertheless, by the result by Weinstein, the only requirement for continuation of periodic orbits from the periodic manifold  $\overline{\Lambda}_{k}$  is that the perturbation be defined in a neighborhood of  $\overline{\Lambda}_{k}$ and is sufficiently small. The computations in \cite{BOZ} therefore remain valid when $\mathcal{U} \subsetneq  \R^{3}$ as long as $E_{k} \subset \mathcal{U}$ holds.

The set $\mathcal{U}$ can be otherwise chosen arbitrarily. For instance we may choose $\mathcal{U}$ small. The above-mentioned inclusion might not hold for all $k \ge 1$. Nevertheless we now remark that this inclusion holds for $k$ sufficiently large.

From $X \in K_{0}^{-1} (0)$ we deduce that 
$$\tau_{k} |z|^{2}+\dfrac{|w|^{2}}{8}=1,$$
which implies that $|z| \le \tau_{k}^{-1/2} \to 0$ as $k \to + \infty$. 

The Hamiltonian function
$$\bar{K}_{\varepsilon}(\bar{X})=K_{\varepsilon}(X)$$
is not well-defined on the whole manifold $\mathcal{M}_{0}$ but it is well-defined on the open set
$$\mathcal{D}_{\delta}=\{ \bar{X}\in \mathcal{M}_0 :\; X=(z,w,\bar{t},\tau),\; |z| < \delta \},$$
where $\delta >0$ is such that
$$\{u \in \R^{3}: |u| < \delta^{2}\} \subset \mathcal{U}.$$

The periodic manifold $\bar{\Lambda}_{k}$ is thus contained in $\mathcal{D}_{\delta}$ for $k$ large enough. For these values of $k$, the argument in \cite{BOZ} holds without change.

We sum up the previous discussions in the following theorem:
\begin{theo} \label{theo: local} Assume that $U: \R \times \mathcal{U} \times [0, \varepsilon_{*}] \to \R$ is a $C^{\infty}$-function satisfying
$$U(t+T, u, \varepsilon)=U(t, u, \varepsilon)$$
for some fixed $T >0$ and an open set $\mathcal{U} \subset \R^{3}$ containing $u=0$. Given an integer $l \ge 1$ there exists $\varepsilon_{l}>0$ such that \eqref{eq:1} has at least $l$ generalized $T$-periodic solutions (lying in $\mathcal{U}$) for each $\varepsilon \in ]0, \varepsilon_{l}[$.
\end{theo}

\begin{rem} The same result holds in the plane and in a line as well, which can be obtained by following the same line of argument and using Levi-Civita regularization instead. 
\end{rem}

In the previous proof we have used that the solutions $X(s)=(z(s), w(s), t(s), \tau(s))$ of the Hamiltonian system $(M, \omega, K_{\varepsilon})$ lying in 
$$\{K_{\varepsilon}(X)=0, BL(X)=0\}$$ 
and satisfying 
$$X(s+S)=g \ast X(s), t(s+S)=t(s)+\eta T, s \in \R$$
for some $S>0, g \in S^{1}$ and $\eta \in \Z$ lead to generalized periodic solutions of \eqref{eq:1} with period $\eta T$. These are given by 
$$u(t)=\overline{z(s(t))} i z(s(t)),$$
where $s=s(t)$ is the inverse of the homeomorphism $t=t(s)$. Note that in the previous proof $\eta=1$ for $\varepsilon=0$ and so $\eta=1$ also for small $\varepsilon$, since $\eta=\eta_{\varepsilon}$ is a continuous function taking values in $\Z$.

The rest of the Section will be devoted to show that it is possible to go from  generalized solutions of \eqref{eq:1} to solutions $X(s)$ of $(M,\omega, K_{\varepsilon})$ satisfying the above conditions. 

A similar discussion for planar solutions and Levi-Civita regularization can be found in Section 4 of \cite{BOZ}. We now explain that the same holds for the spatial case as well. 

Now assume that $u(t)$ is a generalized $T$-periodic solution of  \eqref{eq:1} and consider the function $\sigma(t)=\dfrac{u(t)}{|u(t)|}.$ In principle $\sigma(t)$ is only defined for $t \in \R \setminus \mathcal{Z}$ but the notion of generalized solution  implies that it has a continuous extension to the whole real line. At collisions, the function $\sigma$ is not necessarily smooth but there is some control on its velocity. More precisely for each $t_{0} \in \mathcal{Z}$, there holds
$$\dot{\sigma}(t)=O((t-t_{0})^{-1/3})\,\, \quad \hbox{ as } t \to t_{0}.$$
This asymptotic expansion follows from the formula
$$\dot{\sigma}(t)=\dfrac{|u(t)|^{2} \dot{u}(t)  - \langle u(t), \dot{u}(t) \rangle u(t)}{|u(t)|^{3}}, \qquad t \in \R \setminus \mathcal{Z}$$
together with the classical estimates at collisions (see \cite{Sperling})
$$u(t)=a (t-t_{0})^{2/3} + b(t) (t-t_{0})^{4/3},$$
$$\dot{u}(t)=\dfrac{2}{3}a (t-t_{0})^{-1/3} + c(t) (t-t_{0})^{1/3},$$
where $a \in \R^{3} \setminus \{0\}$ and $b(t), c(t)$ are bounded functions defined in a neighborhood of $t-t_{0}$.
We now identify $\R^{3}$ and $\mathbb{IH}$, $u=u(t), u=u_{1} i + u_{2} j + u_{3} k$ and now view $S^{2}$ and $S^{3}$ as unit spheres in $\mathbb{IH}$ and $\mathbb{H}$ respectively. The properties of $\sigma: \R \to S^{2}$ suggests the following definition.

\begin{defi} Given an interval $[0, T]$ and a partition of it 
$$\mathcal{P}:b_{0}=0 < t_{1 } < \cdots < t_{N}=T,$$
{the path $\alpha:[0,  T] \to S^{d}$, $d=2$ or $3$,} is a $\mathcal{P}$-path, if it is continuous, the restriction $\alpha |_{]t_{i}, t_{i+1}[}$ is $C^{\infty}$ and the integral below is finite:
$$\int_{t_{i}}^{t_{i+1}} |\dot{\alpha} (t) | d t < \infty, i=0, \cdots, N-1.$$
\end{defi}
The map $\Pi :z \mapsto \bar{z} i z$ sends $S^{3}$ onto $S^{2}$ and we can define the Hopf map as the restriction 
$$\pi: S^{3} \to S^{2}, z \mapsto \bar{z} i z.$$
Note that, in the notation of \cite[pp. 376]{Cushman}, $x_{1}=-z_{0}, x_{2}=z_{1}, x_{3}=z_{2}, x_{4}=z_{3}$.

We will need the following result on the lifting of $\mathcal{P}$-paths.
\begin{lem} \label{lem: ehresmann} Assume that $\gamma: [0, T] \to S^{2}$ is a $\mathcal{P}$-path. Then there exists a $\mathcal{P}$-path $\Gamma: [0, T] \to S^{3}$ satisfying $\pi \circ \Gamma=\gamma$ and
$$\Re\{\bar{\Gamma}(t) i \dot{\Gamma}(t)\}=0, \hbox{  when } t \neq t_{i}.$$
Moreover, if $\gamma(0)=\gamma(T)$ then $\Gamma(T)=g \Gamma(0)$ for some $g \in S^{1}$.
\end{lem}

The proof of this result is postponed to the next Section. It is worth to observe that the last conclusion on closed paths follows easily from the  structure of the fibers of the Hopf map
$$\pi^{-1}(\gamma(0))=\{g \Gamma(0): g \in S^{1}\}.$$

Let us take our generalized $T$-periodic solution $u(t)$. The set $[0, T] \cap \mathcal{Z}$ defines a partition $\mathcal{P}$. Then $\sigma$ is a $\mathcal{P}$-path and Lemma \ref{lem: ehresmann} can be applied to $\gamma=\sigma$. In consequence there exists a $\mathcal{P}$-path $\Sigma: [0, T] \to S^{3}$ with $\pi \circ \Sigma=\sigma$ and $\Re\{\bar{\Sigma} i \Sigma\}=\Re\{\bar{
\Sigma} i \dot{\Sigma}\}=0$. Note that $\bar{\Sigma} i \Sigma \in S^{2} \subset \mathbb{IH}$.

Since $\sigma$ is $T$-periodic we extend $\Sigma$ to the whole real line via the formula
$\Sigma(t+T)=g \Sigma(t), t \in \R.$
Then $\Sigma: \R \to S^{3}$ is continuous.

We will define $X(s)=(z(s), w(s), t(s), \tau(s))$ in terms of $\Sigma$. The coordinate $t(s)$  is defined from Sundman integral as in \cite{BOZ}. The set $\mathcal{Z}^{*}=t^{-1} (\mathcal{Z})$ is discrete. Then we have
$$z(s)=r(s) \Sigma(t(s)), s \in \R,$$
where 
$$r(s)=|u(t(s))|^{1/2}, w(s)=4 z'(s), \tau(s)=-E(t(s))+\varepsilon  U(t(s), \bar{z}(s) i z(s), \varepsilon), s \in \R \setminus \mathcal{Z}^{*}.$$ 

After this definition the proof follows along the same line of \cite{BOZ}. The only essential difference is the verification of the additional condition $BL(X(s))=0$ which is equivalent to 
$$\Re\{\bar{z}(s) i z'(s)\}=0.$$
To check this condition, it is sufficient to observe that
$$z'(s)=r'(s) \Sigma(t(s))+r(s)^{3} \dot{\Sigma} (t(s)) $$
and thus the condition follows directly.

\section{Lifting of paths via the Hopf map}\label{33}
For each $z \in S^{3}$ we consider the tangent space
$$T_{z}(S^{3})=\{w \in \mathbb{H}: w \perp z\}$$
with vertical and horizontal subspaces
$$\hbox{Vert}_{z}=\ker (d \Pi )_{z} \cap T_{z} (S^{3})$$
$$\hbox{Hor}_{z}=[\ker (d   \Pi  )_{z}]^{\perp} \cap T_{z} (S^{3})$$
The real vector space $\hbox{Vert}_{z}$ has dimension one and is spanned  by $i z$. The space $\hbox{Hor}_{z}$ has dimension two and is spanned by  the vectors $jz, kz $. Moreover
$$(d \pi)_{z} (\hbox{Hor}_{z})=T_{\pi (z)}(S^{2}).$$
The splitting
$$T_{z}(S^{3})=\hbox{Vert}_{z} \oplus \hbox{Hor}_{z}$$
is clearly smooth with respect to the base points and thus defines an Ehresmann connection associated to the submersion $\pi$.  We refer to \cite[Chapter VIII]{Cushman} for more details.

The results in \cite{Cushman} imply that this connection is ``good''. This means that, given a $C^{\infty}$ curve $\alpha: [t_{0}, t_{1}]  \to S^{2}$ and a point $\xi \in \pi^{-1}  (\alpha (t_{0}))$, there exists a horizontal lift starting at $\xi$. This lift is a $C^{\infty}$ curve $A:[t_{0}, t_{1}] \to S^{3}$ satisfying $A(t_{0})=\xi, \pi \circ A = \alpha, \dot{A}(t) \in \hbox{Hor}_{A(t)}$ for each $t \in [t_{0}, t_{1}]$. 

From our point of view, the key point is the following characterization of the horizontal component
$$\hbox{Hor}_{z}=\{w \in \mathbb{H}: w \perp z, \Re\{\bar{z} i w\}=0\}.$$
From this observation we could derive a version of Lemma \ref{lem: ehresmann} in the class of $C^{\infty}$-paths, or even in the class of piecewise $C^{\infty}$-paths. This second case will follow by an iterative application of the above-mentioned lifting principle for good Ehresmann connections. However we must work in the larger class of $\mathcal{P}$-paths and the proof uses specific properties of the Hopf map. To prove Lemma \ref{lem: ehresmann} we will restrict to the case $\mathcal{P}: t_{0}=0<t_{1}=T$ so that $[t_{0}, t_{1}]=[0, T].$ The iterative argument for a general partition is left to the reader. We also assume that $\gamma (t)\neq i$ for each $t\in [t_0,t_1]$ and, at the end of the proof, we will show that there is no loss of generality.

Following \cite{Cushman}, we shall obtain $\Gamma(t)$ from the solution of a system of differential equations in the plane. 
We sketch some computations which are similar to those in \cite{Cushman}, pp. 376-377. Note that our map $\pi$ does not coincide with the map $\pi$ in \cite{Cushman}, although they are conjugate. 
\par The condition $\dot{\Gamma}(t) \in \hbox{Hor}_{\Gamma(t)}$ implies that
$$\dot{\Gamma}(t)=(\lambda_{1}(t)j + \lambda_{2}(t)k) \Gamma(t) ,$$
where $\lambda_{1}, \lambda_{2}: [t_{0}, t_{1}] \to \R$ are functions to be determined. With the notation $\Gamma=\Gamma_{0} + \Gamma_{1} i + \Gamma_{2} j + \Gamma_{3} k$ we are led to the equations

\begin{equation}\label{eq: one}
\dot{\Gamma}_{0}=-\lambda_{1} \Gamma_{2} -\lambda_{2} \Gamma_{3}, \qquad \dot{\Gamma}_{1}=\lambda_{1} \Gamma_{3} -\lambda_{2} \Gamma_{2}
\end{equation}
\begin{equation}\label{eq: two}
\dot{\Gamma}_{2}=\lambda_{1} \Gamma_{0} +\lambda_{2} \Gamma_{1}, \qquad \dot{\Gamma}_{3}=-\lambda_{1} \Gamma_{1} +\lambda_{2} \Gamma_{0}.
\end{equation}
From the notation $\gamma=\gamma_{1} i +\gamma_{2} j + \gamma_{3} k$ and the definition of $\pi$, 
\begin{equation}\label{eq: arbol}
\gamma_{1}=\Gamma_{0}^{2}+\Gamma_{1}^{2}-\Gamma_{2}^{2}-\Gamma_{3}^{2}, \qquad \begin{pmatrix} \Gamma_{2} & - \Gamma_{3} \\ \Gamma_{3} & \Gamma_{2} \end{pmatrix} \begin{pmatrix} \Gamma_{1}  \\ \Gamma_{0}\end{pmatrix}=\dfrac{1}{2} \begin{pmatrix} \gamma_{2}  \\ \gamma_{3}\end{pmatrix}.
\end{equation}
In particular, $\Gamma_{2}^{2} + \Gamma_{3}^{2}=\dfrac{1}{2} (1-\gamma_{1})$. 

We recall that the curve $\gamma$ does not pass through the point $i$ and therefore $\gamma_1 <1$ and  $\Gamma_{2}^{2}+\Gamma_{3}^{2}>0$ everywhere. This allows us to use the equations in \eqref{eq: one} to solve $\lambda_{1}$ and $\lambda_{2}$. 

Plugging the corresponding formulas in \eqref{eq: two} and taking into account the equations in \eqref{eq: arbol} we obtain
\begin{equation}\label{eq: liason}
(1-\gamma_1 )\begin{pmatrix} \dot{\Gamma}_{2}  \\ \dot{\Gamma}_{3}\end{pmatrix}=-M\begin{pmatrix} \dot{\Gamma}_{1}  \\ \dot{\Gamma}_{0}\end{pmatrix},\; \; M=\begin{pmatrix} \gamma_{2}&\gamma_3 \\ \gamma_3 &-\gamma_2 \end{pmatrix}.
\end{equation}
Also, from the second identity in \eqref{eq: arbol},
\begin{equation}\label{eq:deuxl}
(1-\gamma_1 )\begin{pmatrix} \Gamma_{1}  \\ \Gamma_{0}\end{pmatrix}=  \begin{pmatrix} \Gamma_{2}&\Gamma_3 \\ -\Gamma_3 &\Gamma_2 \end{pmatrix}
\begin{pmatrix}\gamma_2  \\ \gamma_3 \end{pmatrix} =M\begin{pmatrix} \Gamma_{2} \\ \Gamma_3  \end{pmatrix}.
\end{equation}
Differentiating with respect to $t$ and substituting the result in \eqref{eq: liason} we find
$$\begin{pmatrix} \dot{\Gamma}_{2}  \\ \dot{\Gamma}_{3}\end{pmatrix}=-\frac{1}{2(1-\gamma_1)} [\dot{\gamma}_1 M\begin{pmatrix} \Gamma_1  \\ \Gamma_0 \end{pmatrix}+M\dot{M} \begin{pmatrix} \Gamma_{2}  \\ \Gamma_{3}\end{pmatrix}].$$
In these computation we have used that $M^2 =(1-\gamma_1^2)I$. Combining this identity with \eqref{eq:deuxl}, we are led to a planar linear system of the type
\begin{equation}\label{eq: ps}
\begin{pmatrix} \dot{\Gamma}_{2}  \\ \dot{\Gamma}_{3}\end{pmatrix}=B(t) \begin{pmatrix} \Gamma_{2}  \\ \Gamma_{3}\end{pmatrix}.
\end{equation}
where the coefficients of the matrix $B(t)$ are linear combinations of functions of the type $\dfrac{\gamma_{i} \dot{\gamma}_{j}}{1-\gamma_{1}}$ and $ \dfrac{\dot{\gamma}_{i}}{1-\gamma_{1}}$. 

Most probably this matrix is not continuous at $t=t_0=0$ and $t=t_1=T$. 
Nevertheless, since we know that $\gamma $ is $C^{\infty}$ in $]0, T[$ and $\int_{0}^{T} |\dot{\gamma}(t)| d t < \infty$, we deduce that the coefficients of $B(t)$ belong to the Lebesgue space $L^{1}(]0,T[)$ .  We are assuming that $\gamma_{1}(t) \neq 1$ if $t \in [0, T].$ In consequence, the matrix $B$ is integrable and the system \eqref{eq: ps} satisfies the conditions of Carath\'eodory's theorem (see for instance \cite[Chapter 2]{CL}). Given $\xi \in \pi^{-1}(\gamma(0))$, we impose the initial condition 
$$\Gamma_{2}(0)=\xi_{2}, \Gamma_{3}(0)=\xi_{3}$$
to obtain a unique solution of the Cauchy problem for \eqref{eq: ps} defined on the whole interval $[0,T]$. These functions $\Gamma_{2}$ and $\Gamma_{3}$ are absolutely continuous in $[0, T]$. The functions $\Gamma_0$ and $\Gamma_1$ are
defined from the identity \eqref{eq:deuxl} and they are also absolutely continuous. In particular there holds 
$$\int_{0}^{T} |\dot{\Gamma}(t)| d t < \infty$$
and $\Gamma$ is a $\mathcal{P}$-path. Going back to the previous construction we observe that $\Gamma$ is the desired lift. 
\par
To complete the proof we must remove the extra assumption $\gamma_1 \neq 1$. Assume now that $\gamma (t)$ is an arbitrary $\mathcal{P}$-path. Since $\gamma$ is smooth on $]t_0,t_1[$, the set $\gamma ([t_0 ,t_1])$ has zero measure in $S^2$.
Let us take a point $\xi \in S^2$ such that $\gamma (t)\neq \xi$ for each $t\in [t_0 ,t_1]$. We select a rotation of $S^2$ sending $\xi$ into $i$. This rotation can be expressed in the form $z'=qz\overline{q}$ for some $q\in S^3$. Then $\gamma_* (t)
=q\gamma (t)\overline{q}$ is a $\mathcal{P}$-path with $\gamma_* (t)\neq i$ for each $t\in [t_0 ,t_1]$. The possible lifts of $\gamma$ and $\gamma_*$ are linked, for if $\Gamma =\Gamma (t)$ satisfies $\pi \circ \Gamma =\gamma$ then
$\pi \circ \Gamma_* =\gamma_*$, where $\Gamma_* (t)=\Gamma (t)\overline{q}$. Moreover, it is easy to check that $\Re\{\bar{\Gamma}(t) i \dot{\Gamma}(t)\}=0$ is equivalent to $\Re\{\bar{\Gamma}_*(t) i \dot{\Gamma}_*(t)\}=0.$ 

\section{A restricted three-body problem}\label{44}
Let us assume that the $C^{\infty}$ functions $X, x: \R \to \R^{3},  X=X(t), x=x(t),$ are $T$-periodic and satisfy 
$$X(t) \neq x(t)\qquad \forall t \in \R.$$
In addition we assign positive masses $M$ and $m$ to them, {by normalization}
$$M+m=1.$$
We consider the system
\begin{equation}\label{eq: 3R}
\ddot{\xi}=\dfrac{M(X(t)-\xi)}{|X(t)-\xi|^{3}}  + \dfrac{m(x(t)-\xi)}{|x(t)-\xi|^{3}}. 
\end{equation}
This system describes the motion of an infinitesimal body attracted by two moving centers $X(t)$ and $x(t)$. When $(X(t), x(t))$ solves the corresponding two-body problem, then the system describes a restricted spatial three-body problem. By assuming that $(X(t), x(t))$ move on a circular or elliptic Keplerian orbit we obtain a periodic system. Note that for {the} restricted circular three-body problem many studies have been made in a proper rotating coordinate system and for the elliptic problem in a proper {rotating-pulsating coordinate system}. 
We shall just consider the problems in the inertial system. 

In principle we could have collisions of the infinitesimal particle $\xi$  with any of the two primaries: $X=\xi$ or $x=\xi$. Nevertheless we shall only consider collisions with the first primary $X$.  

 A continuous and $T$-periodic function $\xi : \R \to \R^{3},\; \xi=\xi(t)$ is called a generalized periodic solution of the first kind if the following conditions hold:
 \begin{itemize}
 \item $\mathcal{Z}_{M}=\{t \in \R: \xi(t)=X(t)\}$ is discrete;
 \item $\mathcal{Z}_{m}=\{t \in \R: \xi(t)=x(t)\}$ is empty;
 \item In each interval $I \subset \R \setminus \mathcal{Z}_{M}$, the function $\xi(t)$ is $C^{\infty}$ and satisfies \eqref{eq: 3R};
 \item For each $t_{0} \in \mathcal{Z}_{M}$ the limits below exist
 $$\lim_{t  \to t_{0}} \dfrac{\xi(t)-X(t)}{|\xi(t)-X(t)|}, \quad \lim_{t  \to t_{0}} \{\dfrac{1}{2} |\dot{\xi}(t)-\dot{X}(t)|^{2}-\dfrac{M}{|\xi(t)-X(t)|}\}.$$
 \end{itemize}
%Let us now assume that the primaries depend upon a parameter $\varepsilon \in [0, 1]$ so that $X, x: \R \times [0, 1] \to \R^{3}, X=X_{\varepsilon}(t), x=x_{\varepsilon}(t)$ are $C^{\infty}$, $T_{\varepsilon}$-periodic in $t$ and for $t \in \R, \varepsilon \in [0, 1]$ there %holds $X_{\varepsilon}(t) \neq x_{\varepsilon}(t)$.

Let us now assume that the primaries depend upon a parameter $\varepsilon \in [0, 1]$, $X_{\varepsilon}=X_{\varepsilon}(t), x_{\varepsilon}=x_{\varepsilon}(t)$. We assume that $(X_{\varepsilon}(t), x_{\varepsilon}(t))$ is a circular or elliptic solution of the two body problem with masses $M_{\varepsilon}, m_{\varepsilon}$ and with center of mass placed at the origin

\begin{equation}\label{eq:cm}
M_{\varepsilon} X_{\varepsilon} + m_{\varepsilon} x_{\varepsilon}=0. 
\end{equation}

From the general theory of Kepler problem we know that $X_{\varepsilon}(t)$ satisfies
\begin{eqnarray*}
X_{\varepsilon}(t)= R_{\varepsilon} \left(a_{\varepsilon} (\cos u(t)-e_{\varepsilon}), a_{\varepsilon} \sqrt{1-e_{\varepsilon}^{2}} \sin u(t), 0 \right)^* \\
u(t)-e_{\varepsilon} \sin u(t)=\dfrac{m_{\varepsilon}^{3/2}}{a_{\varepsilon}^{3/2}} t.
\end{eqnarray*}
We are assuming that $t=0$ is the time of passage through the pericenter and the matrix $R_{\varepsilon}$ is in the group of rotations $SO(3)$. From $X_{\varepsilon}(t)$ one determines $x_{\varepsilon}(t)$ from the center of mass condition Eq. \eqref{eq:cm}. 

From now on it will be assumed that the functions
$$\varepsilon \in [0, 1] \mapsto R_{\varepsilon} \in SO(3)$$
and
$$\varepsilon \in [0, 1] \mapsto m_{\varepsilon}, M_{\varepsilon}, a_{\varepsilon}, e_{\varepsilon}$$
are all $C^{\infty}$ and, for each $\varepsilon>0$ there holds
$$m_{\varepsilon}, M_{\varepsilon}>0, m_{\varepsilon}+M_{\varepsilon}=1, a_{\varepsilon}>0, \,\, 0 \le e_{\varepsilon} <1.$$

According to the third Kepler law, the system \eqref{eq: 3R} is periodic in time with period

\begin{equation}\label{eq:per}
T_{\varepsilon}=\dfrac{2 \pi}{m_{\varepsilon}^{3/2}} a_{\varepsilon}^{3/2}.
\end{equation}

In the following result we will assume that the primary $X_{0}$ is fixed at the origin while $x_{0}$ describes a circular or elliptic Keplerian orbit with mass $m_{0}=0$. 

\begin{theo}\label{theo: g3R}
Assume in addition that 
$$m_{0} = 0 ,  e_0<1 , T_{\varepsilon} \to T_{0} >0 \hbox{ as } \varepsilon \to 0^{+}. $$
Then , for any given integer $l \ge 1$, there exists $\varepsilon_{l}>0$ such that the equation \eqref{eq: 3R} has at least $l$ generalized $T_{\varepsilon}$-periodic solutions of the first kind for $\varepsilon \in ]0, \varepsilon_{l}[$.
\end{theo}

In contrast to many other results \cite{CPS, AL, {PYFN}}, we do not impose any further resonance or non-degeneracy conditions and the dependence with respect to parameters is very general. On the other hand, our solutions are understood in a generalized sense.

\begin{proof} 
%After restricting the size of $\varepsilon$ we can assume that there exists $\delta>0$ such that
%$$|X_{\varepsilon}(t)-x_{\varepsilon}(t)| \ge \delta \hbox{ for } t \in \R, \varepsilon  \in [0, \varepsilon_{*}].$$
To fix a uniform period  we change time $s \to t$ according to the relation
$$T_{\varepsilon} s=t$$
and set $\eta(s)=\xi(t)$,
so that \eqref{eq: 3R} is transformed into
$$\eta^{''}=\dfrac{T_{\varepsilon}^{2} M_{\varepsilon} (\phi_{\varepsilon}(s)-\eta)}{|\phi_{\varepsilon}(s)-\eta|^{3}}+ \dfrac{T_{\varepsilon}^{2} m_{\varepsilon} (\psi_{\varepsilon}(s)-\eta)}{|\psi_{\varepsilon}(s)-\eta|^{3}},$$
where $\phi_{\varepsilon}(s)=X_{\varepsilon} (T_{\varepsilon} s), \psi_{\varepsilon}(s)=x_{\varepsilon} (T_{\varepsilon} s)$.

We then introduce 
$$u={\lambda^{-1}} (\eta-\phi_{\varepsilon}(s))$$
where $\lambda>0$ is a normalization parameter to be adjusted. In this way we obtain an equation of the form of \eqref{eq:1} if $\lambda=\lambda_{\varepsilon}$ is given by
\begin{equation}\label{eq:para}
\lambda_{\varepsilon}^{3} =T_{\varepsilon}^{2} M_{\varepsilon}.
\end{equation}
Namely,
$$u''=-\dfrac{u}{|u|^{3}}+\varepsilon \nabla_{u} U(s, u, \varepsilon)$$
with
$$U(s, u, \varepsilon)=\dfrac{1}{\varepsilon} \lambda_{\varepsilon}^{-3} T_{\varepsilon}^{2} m_{\varepsilon} \dfrac{1}{|\lambda_{\varepsilon}^{-1} (\psi_{\varepsilon} (s)-\phi_{\varepsilon}(s)) - u|} - \dfrac{1}{\varepsilon} \lambda_{\varepsilon}^{-1} \langle \phi''_{\varepsilon}(s), u \rangle.$$
This function belongs to $C^{\infty} (\R \times \mathcal{U} \times [0, \varepsilon_{*}])$ where $\mathcal{U}$ is a neighborhood of $u=0$ and $\varepsilon_{*}>0$ is small enough. 

To prove this, we first observe that
$$a_{\varepsilon} (1-e_{\varepsilon}) \le |\phi_{\varepsilon} (s)| \le a_{\varepsilon} (1+e_{\varepsilon}), s \in \R.$$
Then, in view of \eqref{eq:cm}, for small enough $\varepsilon$ we have
$$|\lambda_{\varepsilon}^{-1} (\psi_{\varepsilon}(s)-\phi_{\varepsilon}(s))| \ge |\lambda_{\varepsilon}^{-1} \phi_{\varepsilon}(s)| \ge  \lambda_{\varepsilon}^{-1} \dfrac{1}{m_{\varepsilon}} a_{\varepsilon} (1-e_{\varepsilon}).$$ 
Using \eqref{eq:per} and \eqref{eq:para} we see that the lower bound converges to $\dfrac{1-e_{0}}{(2 \pi)^{2/3}}$ as $\varepsilon \to 0$. Consequently there exist $\varepsilon_{*}>0$ and $\delta>0$ such that
$$|\lambda_{\varepsilon}^{-1} (\psi_{\varepsilon} (s) -\phi_{\varepsilon}(s)) | \ge \delta \hbox{ for } \varepsilon \in [0, \varepsilon_{*}], s \in \R.$$ 

We define
$$\mathcal{U}=\{u \in \R^{3}: |u| < \delta\}.$$
From the explicit definition of $\phi_{\varepsilon}$ and $\psi_{\varepsilon}$ it is clear that $U$ is smooth on $\R \times \mathcal{U} \times ]0, \varepsilon_{*}]$. It remains to analyze $\varepsilon=0$.

First we observe that 
$$\phi_{\varepsilon} (s)=a_{\varepsilon} R_{\varepsilon} (\cos u-e_{\varepsilon}, \sqrt{1-e_{\varepsilon}^{2}} \sin u, 0)^*, u-e_{\varepsilon}\sin u=2 \pi s$$
and therefore $$\lambda_{\varepsilon}^{-1} (\psi_{\varepsilon} (s) -\phi_{\varepsilon}(s)) =-\lambda_{\varepsilon}^{-1} (\frac{T_{\varepsilon}}{2\pi})^{2/3} R_{\varepsilon} (\cos u-e_{\varepsilon}, 
\sqrt{1-e_{\varepsilon}^{2}} \sin u, 0)^*$$  is $C^{\infty}$ in $\R\times [0,1]$. 
Also, the function $f(\varepsilon)=\dfrac{1}{\varepsilon} \lambda_{\varepsilon}^{-3} T_{\varepsilon}^{2} m_{\varepsilon}=\dfrac{m_{\varepsilon}}{ \varepsilon M_{\varepsilon}}$ belongs to $C^{\infty}[0, \varepsilon_{*}]$ and therefore the first summand in the definition of $U$ is smooth. 

To analyze the second summand we differentiate twice the function $\phi_{\varepsilon}$. 
Then we have
$$\phi''_{\varepsilon}(s)=a_{\varepsilon} \chi(s, \varepsilon)$$ 
with
$$\chi(s, \varepsilon)=u'' R_{\varepsilon} (-\sin u, \sqrt{1-e_{\varepsilon}^{2}} \cos u,0)^*-(u')^{2} R_{\varepsilon}  (\cos u, \sqrt{1-e_{\varepsilon}^{2}} \sin u, 0)^*   $$
and
$$u'=\dfrac{\partial u}{\partial s}=\dfrac{2 \pi}{1-e_{\varepsilon} \cos u}, u''=\dfrac{\partial^{2} u}{\partial s^{2}}=-\dfrac{2 \pi e_{\varepsilon} \sin u}{(1-e_{\varepsilon} \cos u)^{2}} u'. $$
The function $\chi(s, \varepsilon)$ thus belongs to $C^{\infty} (\R \times [0, \varepsilon_{*}], \R^{3})$. In addition, 
$$g(\varepsilon)=\dfrac{1}{\varepsilon} \lambda_{\varepsilon}^{-1} a_{\varepsilon}=\dfrac{m_{\varepsilon}}{\varepsilon M_{\varepsilon}^{1/3} (2 \pi)^{2/3}}$$ is in $C^{\infty}[0, \varepsilon_{*}]$.

Therefore $U$ is smooth in $\R \times \mathcal{U} \times [0, \varepsilon_{*}]$ and Theorem \ref{theo: local} is applicable. Undoing the change of variables we obtain generalized $T_{\varepsilon}-$periodic solutions. To check the continuity of the energy
at collisions it is convenient to use the identity $$\frac{1}{2} |u'(s)|^2 -\frac{1}{|u(s)|} =\frac{T_{\varepsilon}^2}{\lambda_{\varepsilon}^2}\Bigl(\frac{1}{2} |\dot{\xi}(t)-\dot{X}_{\varepsilon} (t)|^2 -\frac{M_{\varepsilon}}{|\xi (t)-X_{\varepsilon}(t)|}\Bigr).$$ 
\end{proof}

\begin{rem} After the change of variables $\xi =R_{\varepsilon} \xi_1$ we can assume that the primaries lie on the fixed plane $x_3=0$. Then we can apply the planar version of  Theorem \ref{theo: local} to conclude that, for each $\varepsilon$, the three bodies move in a common plane. 
\end{rem}

We end this note with a comparison of our result with classical results concerning periodic orbits of the planar circular restricted three-body problem of the first and second kind in a rotating frame, where first and second kind refer to continuations of periodic orbits from circular and elliptic periodic orbits of the limiting Kepler problem in a rotating frame respectively. Indeed the period of rotations of the reference frame is derived from the corresponding period $T$ of the Keplerian elliptic motions of the primaries. When a periodic orbit has minimal period $T/n$ in the rotating frame, then { it will also have }period $T$ and therefore {this period is preserved} in the initial fixed reference frame. The existence of such orbits of the first kind has been obtained by Poincar\'e \cite{P} and has been explained in \cite{MZ}. 
 
 \begin{ack} R. O. is supported by  MTM2017-82348-C2-1-P, L. Z. is supported by DFG ZH 605/1-1. 
 \end{ack}

\medskip
\medskip
\medskip
\medskip
\medskip

\noindent
Rafael Ortega, Departamento de Matem\'atica Aplicada, Universidad de Granada, 18071 Granada, Spain: rortega@ugr.es
\newline
Lei Zhao, Institut f\"ur Mathematik, Universit\"at Augsburg, Universit\"atsstra{\ss}e 2, 86159 Augsburg, Germany: lei.zhao@math.uni-augsburg.de

\end{document}